\documentclass{article}
\usepackage[utf8]{inputenc}
\usepackage{comment}
\usepackage{changepage}
\usepackage[colorlinks,citecolor=blue,urlcolor=blue,bookmarks=false,hypertexnames=true]{hyperref} 
\title{Properties of terms of OEIS A342810}
\author{R\"udiger Jehn and Kester Habermann}
\date{June 2021}

\begin{document}

\maketitle

{\bf Abstract}\\

The OEIS sequence A342810 contains the numbers that divide the smallest number that has the sum of their digits. It is proved that if a term $x$ has the form $3^m \times y$ where $m \geq 2$, then all prime factors of $y$ are prime divisors of solutions to $10^n \equiv 1 (\bmod n)$. It is also proved that if a term $x$ has the form $3^m \times p \times q$ where $m \geq 2$, where $p$ is a prime divisor of a solution to $10^n \equiv 1 (\bmod n)$ and where $q$ is the product of all other factors of the prime factorisation of $x$, then all numbers $3^m \times p^i \times q$ are also terms of the sequence A342810 for any integer~$i$.

\section{Introduction}

The OEIS sequence A342810 contains the numbers that divide the smallest number that has the sum of their digits. For instance 27 is a term of this sequence because 27 divides 999 which is the smallest number with a digit sum of 27.\\

21, 27, 81, 191, 243, 729, 999, 2187, 2997, 6561, 8991, 19683, 26973, 33321, 36963, 39049, 59049 and 80919 are the terms of this sequence between 10 and 100000. There are 811 terms smaller than $10^{13}$.\\

It is straightforward to show that none of the terms larger than 10 of A342810 is divisible by 2, 5, 11 and 13. This is proved in Chapter 2.\\

It can be observed that if a term $x$ has the form $3^m \times y$ where $m \geq 2$ (which is the case for the overwhelming number of terms of this sequence), then all prime factors of $y$ are in the OEIS sequence A066364 \cite{A066364}. This sequence 3, 37, 163, 757, ... contains the prime divisors of solutions to $10^n \equiv 1 (\bmod n)$. This is proved in Chapter 3.\\

It can also be observed that if a term $x$ has the form $3^m \times p \times q$ where $m \geq 2$, where $p$ is a term of OEIS sequence A066364 \cite{A066364} and where $q$ is the product of all other factors of the prime factorisation of $x$, then all numbers $3^m \times p^i \times q$ are also terms of the sequence for any integer $i$. This is proved in Chapter 4.\\

Finally in Chapter 5, the rare cases of terms not divisible by 9 are discussed.

\section{Terms are not divisible by 2, 5, 11 and 13}

If a term $x$ is larger than 10, the sum of its digits is larger than 10 and the smallest number $k$ with a sum of digits of $x$ ends with a "9". If $x$ is a multiple of 2 or 5 it cannot divide a number ending with a "9". Hence no term of sequence A342810 is divisible by 2 or 5.\\

The smallest number $k$ with a sum of digits of $x$ is divisible by 11 only if the alternating digit sum of $k$ is divisible by 11. But there is only one possibility: $k$ consists of an even number of "9"s and the alternating digit sum is 0. The digit sum itself is 9 times an even number and therefore is even and cannot divide the odd number $k$. Hence no term of sequence A342810 is divisible by 11.\\

The smallest number $k$ with a sum of digits of $x$ is divisible by 13 only if the alternating sum of blocks of three digits of $k$ is divisible by 13. There are these four possibilities: 
\begin{enumerate}
    \item $k = (999 999)_n$
    \item $k =  39 (999 999)_n$
    \item $k = 299 (999 999)_n$
    \item $k = 89 999 (999 999)_n$
\end{enumerate}
where $(999 999)_n$ denotes $n$ blocks with 6 "niners" ($n$ being any non-negative integer).
The digit sum $x$ of $k$ is 0, 12, 20 or 44 + $54n$ which is even and cannot divide the odd number $k$. Hence no term of sequence A342810 is divisible by 13.\\

\section{Prime factors of terms divisible by 9}

If $x$ is divisible by 9, then the smallest number that has the sum of their digits is the number $k$ that consists of $\frac{x}{9}$ niners: 9999...
If we write $x = 9pq$ where $p$ is a prime and $q$ can be any composite number including a "1" then $k = 10 ^{pq} - 1$. $x$ divides $k$ and therefore we need to search for primes $p$ for which
$$pq | 10 ^{pq} - 1$$

This means $ 10 ^{pq} \equiv 1 (\bmod pq)$ and we deduce:
$$ 10 ^{pq} \equiv 1 (\bmod pq) \Rightarrow 10 ^{pq} \equiv 1 (\bmod p) $$

Let $e = O_p(10)$ be the multiplicative order of 10 modulo p, which means that the exponent $pq$ must be a multiple of $e$:

$$e | pq $$

If $e$ is even (which is true for $p$ = 7, 11, 13, 17, 19, 23, 29, 47, ...), then an even number $e$ would divide the odd number $10 ^{pq} - 1$, which is a contradiction. Hence primes with an even multiplicative order of 10 cannot be a divisor of a term $x$ that is divisible by 9.\\

In general if $x = 9n$ and $10^n \equiv 1 (\bmod n)$ we know that the prime factors of $n$ are terms of the sequence A066364, because this is the definition of this sequence.\\

Already in 1989 Kennedy and Cooper \cite{kennedy} published a paper analysing solutions of the equation $10^n \equiv 1 (\bmod n)$. The prime factors $p$ of the solutions are listed in Table \ref{table:1}:

\begin{table}[h!]
\centering
\begin{tabular}[h]{|c|c|c|}
\hline
$p$ & $e$ & $n_p$ \\
\hline
3 & 1 & 0 \\
37 & 3 & 1 \\
163 & $81 = 3^4$ & 4 \\
757 & $27 = 3^3$ & 3 \\
1999 & $999 = 3^3 \cdot 37$ & 3 \\
\hline 
\end{tabular}
\caption{Prime factors of terms divisible by 9, their multiplicative order of 10 and the exponent $n_p$ of 3 in the prime factorisation of $e$ which will be used in the next chapter.}
\label{table:1}
\end{table}

\section{Deriving new terms by increasing the exponent of a prime factor}

In the previous chapter we found all prime factors of the terms divisible by 9. For these terms this theorem holds:\\

\begin{adjustwidth}{1cm}{}
If a term $x$ has the form $3^m \times p \times q$ where $m \geq 2$, where $p$ is a term of OEIS sequence A066364 \cite{A066364} and where $q$ is the product of all other factors of the prime factorisation of $x$, then all numbers $3^m \times p^i \times q$ are also terms of the sequence for any integer $i$.\\
\end{adjustwidth}

For the proof, we need the divisibility rules for the numbers of the sequence A066364. For 3 it is well known: a number is divisible by 3 if the digit sum is divisible by 3. A less known rule exists for 37: a number is divisible by 37 if the 3-block digit sum is divisible by 37. (You obtain the 3-block digit sum of a number by splitting the number in blocks of three digits starting from the right and summing up the blocks.) And in general, there are similar rules for all numbers $y$ of the OEIS sequence A066364 since there exist an $e$ for which $y | (10^{e} - 1)$. In fact $e$ is the multiplicative order of 10 modulo y denoted as $O_y(10)$ as listed in Table \ref{table:1}.\\

The divisibility rule is: a number is divisible by $y$ if the $O_y(10)$-block digit sum is divisible by $y$.
For instance $e_{757} = 27$ means that a number is divisible by 757 if the 27-block digit sum is divisible by 757.
\\

Assume $x = 3^{n_p+2} \times p \times q$ is a term of A342810. Then $k_x$, the smallest number that has the sum of digits $x$ is
$$ k_x = "9 \ldots 9" \hspace{6mm}  (3^{n_p} \times p \times q \mbox{  times "9")}$$
and
$$ \frac{k_x}{x} = b \mbox{  with } b \mbox{  integer.}$$

Now we investigate $z = p \times x$. Then $k_z$, the smallest number that has the sum of digits $z$ is
$$ k_z = "k_x k_x \ldots k_x"$$
with the block $k_x$ repeating $p$ times. Since $\frac{k_x}{x} = b$ we obtain
$$ \frac{k_z}{x} = "b b \ldots b"$$
where $b$ might need to be completed with leading zeros to have the same string length as $k_x$. For example, if $p=3$, $x=81$ and $k_x = 999,999,999$ then $b = 012,345,679$.\\

After adding the leading zeros, the length of the string $b$ is the same as the length of $k_x$ which is $3^{n_p} \times p \times q$. $b$ can be split in $p$ blocks of length $3^{n_p} \times q$: 
$$b = b_1 b_2 \ldots b_p$$
Each $b_i$ can be split further into $q$ blocks of length $3^{n_p}$:
$$b_i = b_{i1} b_{i2} \ldots b_{iq}$$
Thus the $3^{n_p}$-block digit sum of $"b b \ldots b"$ becomes: 

$$p \times \Sigma_{i=1}^p \Sigma_{j=1}^q b_{ij} $$

which is obviously divisible by $p$.\\

By the divisibility rule also $\frac{k_z}{x}$ is divisible by $p$ and from this it follows immediately that
$$ p \times x | k_z \mbox{ or } z | k_z$$
which means that $z$ is a term of A342810. \\

This incurs that the sequence A342810 has an infinite number of terms.\\

\section{Terms not divisible by 9}

The first 20 terms that are divisible by 3 but not by 9 are listed in Table \ref{table:2} together with their prime factorisation. There is no obvious patter in the prime factors.\\

\begin{table}[h!]
\centering
\begin{tabular}[h]{|r|l|}
\hline
x & Prime factors \\
\hline
3 & 3 \\
6 & 3 * 2 \\
21 & 3 * 7 \\
33321 & 3 * 29 * 383 \\
100389 & 3 * 109 * 307 \\
177897 & 3 * 19 * 3121 \\
       7887189 & 3 * 43 * 61141 \\
       9972201 & 3 * 29 * 83 * 1381 \\
      42874617 & 3 * 179 * 79841 \\
     596535879 & 3 * 131 * 1049 * 1447 \\
    4386835767 & 3 * 19 * 97 * 193 * 4111 \\
    5987941653 & 3 * 1129 * 1767919 \\
   25584896001 & 3 * 7 * 4229 * 288089 \\
   66729175779 & 3 * 7 * 383 * 8296553 \\
 2300294146809 & 3 * 29 * 26440162607 \\
 2348781662217 & 3 * 106871 * 7325909 \\
 6679356919437 & 3 * 19 * 61 * 1921011481 \\
17706895960503 & 3 * 23 * 727 * 14867 * 23743 \\
18436608407469 & 3 * 7 * 541 * 1622798029 \\
36525783703737 & 3 * 439 * 96179 * 288359 \\
\hline 
\end{tabular}
\caption{The first 20 terms of sequence A342810 that are divisible by 3 but not by 9.}
\label{table:2}
\end{table}

Terms not divisible by 3 are extremely rare. Between 10 and $10^{13}$, we found only 5. They are listed in Table \ref{table:3}. It is remarkable that the remainder of the division of a term by 9 is always 2 or 7, but this maybe just a coincidence. Up to $10^{13}$ no other remainders were found.

\begin{table}[h!]
\centering
\begin{tabular}[h]{|r|l|c|}
\hline
x & Prime factors & x mod 9 \\
\hline
           191 & 191 & 2 \\
         39049 & 17 * 2297 & 7 \\
   13778099993 & 7 * 1968299999 & 2 \\
   64894173577 & 26113 * 2485129 & 7 \\
  196742723591 & 19 * 10354880189 & 2 \\
\hline 
\end{tabular}
\caption{The first 5 terms of sequence A342810 that are larger than 10 and not divisible by 3.}
\label{table:3}
\end{table}
\bibliography{references} 
\bibliographystyle{ieeetr}

\end{document}